\documentclass[10pt,oneside]{amsart}

\usepackage{amsmath}
\usepackage{amssymb}
\usepackage{amscd}

\usepackage{verbatim}

\title{Galois points for a normal hypersurface}
\author{Satoru Fukasawa \& Takeshi Takahashi}

\address{Department of Mathematics,  
School of Science and Engineering, Waseda University, 
Ohkubo 3-4-1, Shinjuku, Tokyo 169-8555, Japan.}
\email{s.fukasawa@kurenai.waseda.jp}
\address{Division of General Education, Nagaoka National College of Technology, 888 Nishikatakai, Nagaoka, Niigata 940-8532, Japan}
\email{takeshi@nagaoka-ct.ac.jp}

\newtheorem{theorem}{Theorem}
\newtheorem{proposition}{Proposition}
\newtheorem{corollary}{Corollary}
\newtheorem{lemma}{Lemma}

\theoremstyle{definition}
\newtheorem{example}{Example}
\newtheorem{remark}{Remark}
\newtheorem{note}{Note} 

\newtheorem{definition}{Definition} 

\newcommand{\Gal}{\mathop{\mathrm{Gal}}\nolimits}

\begin{document}
\begin{abstract} 
We study Galois points for a hypersurface $X$ with $\dim {\rm Sing}(X) \le \dim X-2$. 
The purpose of this article is to determine the set $\Delta(X)$ of Galois points in characteristic zero: Indeed, we give a sharp upper bound of the number of Galois points in terms of $\dim X$ and $\dim {\rm Sing}(X)$ if $\Delta(X)$ is a finite set, and prove that $X$ is a cone if $\Delta(X)$ is infinite.  
To achieve our purpose, we need a certain hyperplane section theorem on Galois point. 
We prove this theorem in arbitrary characteristic. 
On the other hand, the hyperplane section theorem has other important applications: For example, we can classify the Galois group induced from a Galois point in arbitrary characteristic and determine the distribution of Galois points for a Fermat hypersurface of degree $p^e+1$ in characteristic $p>0$. 
\end{abstract}
\maketitle

\section{Introduction}  
Let the base field $K$ be an algebraically closed field of characteristic $p \ge 0$ and let $X \subset \Bbb P^{n+1}$ be an irreducible and reduced hypersurface of dimension $n$ and of degree $d \ge 4$ with $s(X):=\dim {\rm Sing}(X)$. 
H. Yoshihara introduced the notion of {\it Galois point} (see \cite{fukasawa2, miura-yoshihara, yoshihara1, yoshihara2, yoshihara3}). 
If the function field extension $K(X)/K(\Bbb P^n)$, induced from the projection $\pi_P:X \dashrightarrow \Bbb P^n$ from a point $P \in \Bbb P^{n+1}$, is Galois, then the point $P$ is said to be Galois. 
In this paper a {\it Galois point} $P$ means {\it Galois point which is contained in the smooth locus $X_{\rm sm}$ or $\Bbb P^{n+1} \setminus X$}, except for Subsection 3.1. 
Let $\Delta(X)$ (resp. $\Delta'(X)$) be the set of all Galois points contained in $X_{\rm sm}$ (resp. $\Bbb P^{n+1} \setminus X$). 
We are interested in the sets $\Delta(X)$ and $\Delta'(X)$. 
If $p=0$, $n=1$ and $X$ is smooth, then Yoshihara determined the sets $\Delta(X)$ and $\Delta'(X)$ (\cite{yoshihara1}). 
If $p>0$, $n=1$ and $X$ is smooth, then the first author determined in most cases (\cite{fukasawa1, fukasawa2}). 
If $p=0$ and $X$ is a smooth quartic surface, then Yoshihara determined the set $\Delta(X)$ completely (\cite{yoshihara2}). 
If $K=\Bbb C$, $n \ge 2$ and $X$ is smooth, then Yoshihara gave sharp upper bounds on the cardinalities of $\Delta(X)$ and $\Delta'(X)$ (\cite{yoshihara3}). 
If $p=0$ and $X$ is a normal quartic surface, then the second author determined the set $\Delta(X)$ (\cite{takahashi}). 
However, for a higher dimensional hypersurface with large singularities, there are few results. 

The purpose of this article is to determine the sets $\Delta(X)$ and $\Delta'(X)$ for a hypersurface of dimension $n$ with $s(X)=\dim {\rm Sing}(X) \le n-2$ in $p=0$. 
We define the dimension of the empty set as $-1$.  
Our results are: 

\begin{theorem} \label{InnerGalois}
Let $X \subset \Bbb P^{n+1}$ be a hypersurface of dimension $n \ge 1$ and degree $d \ge 4$ with $s=s(X) \le n-2$ in characteristic $p=0$ and let $m=m(n,s):= [(n+s+1)/2]$ where $[*]$ means the integer part of $*$. 
Assume that $\Delta(X)$ is a finite set. 
\begin{itemize}
\item[(I-0)] If $d=4$ and $s=-1$ then $\sharp\Delta(X) \le 4(m(n,s)+1)$. 
The equality holds if and only if $X$ is projectively equivalent to the hypersurface defined by 
$$ X_{m+1}X_0^3+\dots+X_{2m+1}X_{m}^3+X_{m+1}^4+\dots+X_{n+1}^4=0.$$ 
\item[(I-1)] If $d=4$, $s \ge 0$ and $n+s$ is odd, then $\sharp\Delta(X) \le 4(m(n,s)-s-1)+(s+2)$. 
The equality holds if and only if $X$ is projectively equivalent to the hypersurface defined by \\
\ \ $X_{m+1}X_0^3+\dots+X_{2m-s-1}X_{m-s-2}^3+X_{2m-s}(X_{m-s-1}^3+\dots+X_{m}^3)$ \\
\hspace{70mm} $+X_{m+1}^4+\dots+X_{2m-s-1}^4+aX_{2m-s}^4=0 $ \\
where $a=0$ or $1$. 
\item[(I-2)] If $d=4$ and $s=0$ and $n$ is even then $\sharp\Delta(X) \le 4m(n,s)+1$. 
The equality holds if and only if $X$ is projectively equivalent to the hypersurface defined by 
$$ X_{m+1}X_0^3+\dots+X_{2m+1}X_{m}^3+X_{m+1}^4+\dots+X_{2m}^4=0. $$
\item[(I-3)] If $d=4$ and $s=1$ and $n+s$ is even then $\sharp\Delta(X) \le 4(m(n,s)-1)$. 
The equality holds if and only if $X$ is projectively equivalent to the hypersurface defined by 
$$ X_{m-1}X_0^3+\dots+X_{2m-3}X_{m-2}^3+X_{m-1}^4+\dots+X_{2m-3}^4+G=0 $$
where $G \in K[X_{2m-2}, X_{2m-1}, X_{n+1}]$ has a multiple component. 
\item[(I-4)] If $d=4$, $s \ge 2$ and $n+s$ is even, then $\sharp\Delta(X) \le 4(m(n,s)-s-1)+(s+2)$. 
Furthermore, the equality holds if and only if $X$ is projectively equivalent to the hypersurface defined by one of the followings: 
\begin{itemize}
\item[(i)] $X_{m-s}X_0^3+\dots+X_{2m-2s-1}X_{m-s-1}^3+X_{m-s}^4+\dots+X_{2m-2s-1}^4+G_1=0$ \\ 
where $s=2$ and $G_1 \in K[X_{2m-2s}, \ldots, X_{n+1}]$ has a multiple component; or  
\item[(ii)] $ X_{m+1}X_0^3+\dots+X_{2m-s-1}X_{m-s-2}^3+X_{2m-s}X_{m-s-1}^3+A_{m-s}X_{m-s}\dots+A_mX_m^3$ \\
\hspace{75mm} $+X_{m+1}^4+\dots+X_{2m-s-1}^4+G_2 =0 $ \\    
where $s \ge 2$ and $A_{m-s}, \ldots, A_m, G_2 \in K[X_{2m-s}, X_{n+1}]$. 
\end{itemize} 
\item[(II)] If $d \ge 5$, then 
$ \sharp \Delta(X) \le m(n,s)+1$.  
Furthermore, the equality holds if and only if $X$ is projectively equivalent to the hypersurface defined by 
$$ X_{m+1}X_0^{d-1}+\dots+X_{2m-s}X_{m-s-1}^{d-1}+A_{m-s}X_{m-s}^{d-1}+\dots+A_mX_m^{d-1}+G=0 $$
where $A_{m-s}, \ldots, A_m, G \in K[X_{m+1}, \ldots, X_{n+1}]$ are homogeneous polynomials with $\deg A_{m-s}=\dots=\deg A_m=1$ and $\deg G=d$.  
(Note that the polynomial $A_{m-s}X_{m-s}+\dots+A_{m}X_m^{d-1}$ does not appear if $s=-1$.) 
\end{itemize} 
\end{theorem} 

\begin{corollary} 
Let $\lambda:=[n/2]$ and 
let $X \subset \Bbb P^{n+1}$ be a hypersurface of dimension $n \ge 1$ and degree $d = 4$ with $s(X) \le n-2$ in characteristic $p=0$. 
Assume that $\Delta(X)$ is a finite set. 
Then, $\sharp\Delta(X) \le 4(\lambda+1)$. 
Furthermore, the equality holds if and only if $X$ is projectively equivalent to the hypersurface defined by 
$$ X_{\lambda+1}X_0^3+\dots+X_{2\lambda+1}X_{\lambda}^3+X_{\lambda+1}^4+\dots+X_{n+1}^4=0.$$ 
\end{corollary}

\begin{theorem} \label{OuterGalois}
Assume that $p=0$, $s(X) \le n-2$ and $\Delta'(X)$ is finite. 
\begin{itemize}
\item[(1)] If $s=-1$, then $\sharp\Delta'(X) \le n+2$. 
Furthermore, the equality holds if and only if $X$ is projectively equivalent to the Fermat hypersurface. 
\item[(2)] If $s \ge 0$, then $\sharp\Delta'(X) \le n-s$. 
Furthermore, the equality holds if and only if $X$ is projectively equivalent to the hypersurface defined by 
$$ X_0^d+\dots+X_{n-s-1}^d+G=0 $$
where $G \in K[X_{n-s}, \ldots, X_{n+1}]$ has a multiple component.  
\end{itemize} 
\end{theorem}

\begin{theorem} \label{InfiniteGalois}
Assume that $p=0$ and $s(X) \le n-2$. 
If the set $\Delta(X) \cup \Delta'(X)$ is infinite, then $X$ is a cone. 
Furthermore, there exist linear spaces $M_1$, $M_2 \subset \Bbb P^{n+1}$ with $M_1 \cap M_2=\emptyset$ and an irreducible hypersurface $Y \subset M_1$ satisfying the following conditions: 
$(i)$ $X=Y \sharp M_2$;
$(ii)$ $M_2$ has the maximal dimension as a vertex of $X$;
$(iii)$ the set $\Delta(Y) \cup \Delta'(Y)$ is finite; and 
$$ (iv) \ \Delta(X)=(\Delta(Y) \sharp M_2) \setminus M_2 \mbox{ and } \Delta'(X)=(\Delta'(Y) \sharp M_2) \setminus M_2, $$
where $\sharp$ means a linear join of algebraic sets. 
\end{theorem}

To prove our Theorems, we need a hyperplane section theorem below. 
We denote by $G_P(X)$ the Galois group induced from $P$ if $P$ is a Galois point for a hypersurface $X$.  

\begin{theorem} \label{HyperplaneSection}
Let $X \subset \Bbb P^{n+1}$ ($n \ge 2$) be an irreducible hypersurface of degree $d \ge 4$ with $0 \le s(X) \le n-2$ (resp. $s(X)=-1$) in characteristic $p \ge 0$. 
Let $P \in \Bbb P^{n+1}$ be a Galois point for $X$. 
Then, the followings hold: 
\begin{itemize}
\item[(i)] A general hyperplane $H$ passing through $P$ satisfies the following condition: 
\begin{itemize}
\item [$(\star)$] the hyperplane section $X_H:=X \cap H$ is an irreducible hypersurface in $H \cong \Bbb P^n$ of degree $d$ with $s(X_H)$ $=s(X)-1$ (resp. $s(X_H)=-1$), and a general tangent space of $X_H$ does not contain $P$.  
\end{itemize} 
\item[(ii)] Let $H$ be a hyperplane passing through $P$ and satisfying the condition $(\star)$. 
Then, the point $P$ is Galois with respect to $X_H$.  
\item[(iii)] The Galois groups are isomorphic: $G_P(X) \cong G_P(X_H)$ for any hyperplane $H$ passing through $P$ and satisfying $(\star)$. 
\end{itemize}   
\end{theorem} 

In fact, this Theorem will be used for the proof of Propositions \ref{l2} and \ref{l2-2} in Subsection 3.1. 
It is remarkable, as in these Propositions, that the birational transformation induced by $\sigma \in G_P(X)$ can be extended to a projective transformation on $\Bbb P^{n+1}$ when $p=0$.  

The hyperplane section theorem has other important applications. 
One application is to classify the Galois group. 
Combining Theorem \ref{HyperplaneSection} with a result on the structure of $G_P(C)$ for a smooth plane curve $C$ with a Galois point $P$ obtained in \cite{fukasawa1}, we have the following: 

\begin{corollary} \label{GroupStructure}
Let $X \subset \Bbb P^{n+1}$ be a hypersurface of degree $d \ge 4$ with $s(X) \le n-2$ in characteristic $p \ge 0$. 
Let $P \in X$ (resp. $P \in \Bbb P^{n+1} \setminus X$) be a Galois point for $X$ and $d-1=p^el$ (resp. $d=p^el$), where $l$ is not divisible by $p$ if $p>0$.  

Then, $l$ divides $p^e-1$ and the Galois group $G_P$ is isomorphic to the semidirect product of $(\Bbb Z/p\Bbb Z)^{\oplus e}$ (as a normal subgroup) and $\Bbb Z/l \Bbb Z$ (as a quotient). 
\end{corollary}

Another application is to determine the distribution of Galois points for a Fermat hypersurface $F_n(q+1) \subset \Bbb P^{n+1}$ of degree $q+1 \ge 4$ in $p>0$ where $q$ is a power of $p$. 
Homma \cite{homma} determined the distribution of Galois points for a Hermitian curve in $\Bbb P^2$, which is defined by $X^{q}Z+XZ^{q}-Y^{q+1}=0$. 
Note that the Hermitian curve is projectively equivalent to $F_1(q+1)$ over $\Bbb F_{q^2}$. 
We will generalize Homma's result. 

\begin{theorem} \label{Fermat} 
A point $P \in \Bbb P^{n+1}$ is Galois for $F_n(q+1)$ if and only if $P \in \Bbb P^{n+1}$ is $\Bbb F_{q^2}$-rational. 
\end{theorem} 

Note that Theorem \ref{Fermat} implies that Theorems \ref{InnerGalois} and \ref{OuterGalois} do not hold in $p>0$. 
In the final section, we give examples which are not cones and have infinitely many Galois points. 
Especially, Example \ref{example1} implies that Theorem \ref{InfiniteGalois} does not hold in $p>0$.

\section{Proof of Theorem \ref{HyperplaneSection}} 
In this section, $X \subset \Bbb P^{n+1}$ is an irreducible hypersurface of degree $d \ge 4$ with $s(X)=\dim {\rm Sing}(X) \le n-2$, otherwise specified. 
We denote by $\check{\Bbb P}^{n+1}$ the dual projective space, which parameterizes hyperplanes in $\Bbb P^{n+1}$. 
We denote by $\Bbb T_xX \subset \Bbb P^{n+1}$ the projective tangent space at $x \in X_{\rm sm}$. 
We have the Gauss map $\gamma: X_{\rm sm} \rightarrow \check{\Bbb P}^{n+1}$; $x \mapsto \Bbb T_xX$. 
We note the following elementary fact: 

\begin{note} 
A hyperplane $H$ coincides with the tangent space $\Bbb T_xX \subset \Bbb P^{n+1}$ at a smooth point $x \in X$ if and only if the scheme $X \cap H$ is singular at $x \in X \cap H$. 
\end{note} 

Let $V_P \subset \check{\Bbb P}^{n+1}$ be the set of all hyperplanes which pass through $P$ and let $S_X \subset \check{\Bbb P}^{n+1}$ be the set of all hyperplanes which contain some irreducible component $Y$ of the singular locus of $X$ with $\dim Y=s(X)$.  
We have the following assertion of Bertini type in arbitrary characteristic (cf. the proof of \cite[II, 8.18]{hartshorne}). 

\begin{note}
Assume that $0 \le s(X) \le n-2$ (resp. $s(X)=-1$). 
For any hyperplane $H \in V_P \setminus (\gamma(X_{\rm sm}) \cup S_X)$, 
the hyperplane cut $X_H:=X \cap H$ is an irreducible hypersurface of degree $d$ with $s(X_H)=s(X)-1$ (resp. $s(X_H)=-1$). 
\end{note} 

A point $P \in \Bbb P^{n+1}$ is called a strange center if the tangent space $\Bbb T_xX$ contain $P$ for a general point $x \in X$. 

\begin{note}
The projection $\pi_P:X \dashrightarrow \Bbb P^n$ is generically finite and separable if and only if $P \in \Bbb P^{n+1}$ is not a strange center. 
\end{note} 

Now we assume that $\pi_P$ is separable and generically finite onto its image. 
Then, the differential map $d_{Q}\pi_P:T_{Q}X \rightarrow T_{\pi_P(Q)}\Bbb P^n$ is isomorphic at a general point $Q \in X$, where $T_QX$ is the Zariski tangent space at $Q$. 
We denote by $U_P$ the maximal open set of $X_{\rm sm}$ such that the differential map $d_{Q}\pi_P:T_{Q}X \rightarrow T_{\pi_P(Q)}\Bbb P^n$ is isomorphic for any $Q \in U_P$. 
We denote by $\Sigma_P \subset \check{\Bbb P}^{n+1}$ the finite set consisting of hyperplanes $H$ such that $X \cap H$ is contained in $X \setminus U_P$ (as a set). 
On the other hand, we note the following: 
\begin{note} 
If $X \cap H$ is an integral scheme, then $\Bbb T_Q(X_H)=\Bbb T_QX \cap H$ for any smooth point $Q$ of $X_H$. 
\end{note} 

Combining above Notes, we have the following: 

\begin{lemma}
Let $X \subset \Bbb P^{n+1}$ be an irreducible hypersurface of degree $d$ with $0 \le s(X) \le n-2$ (resp. $s(X)=-1$). 
Assume that $P$ is not a strange center. 
Then, we have the followings: 
\begin{itemize}
\item[(i)] $V_P \setminus (\gamma(X_{\rm sm}) \cup S_X \cup \Sigma_P) \ne \emptyset$. 
\item[(ii)] Let $H \in V_P \setminus (\gamma(X_{\rm sm}) \cup S_X \cup \Sigma_P)$.  
Then, $H$ satisfies the condition $(\star)$.  
\end{itemize} 
\end{lemma}

Now assume that $P$ is a Galois point.  
Let $G_P$ be the Galois group. 
Then, we find easily that $P$ is not a strange center. 
Therefore, $V_P \setminus (\gamma(X_{\rm sm}) \cup S_X \cup \Sigma_P) \ne \emptyset$ as in Lemma 1. 
Then, we can consider $\sigma \in G_P$ as a birational map from $X$ to itself. 
Needless to say, $\sigma$ is defined and isomorphic at a general point of $X$. 
We denote by the maximal open subset $U_{\sigma}$ such that $\sigma$ is defined and isomorphic over $U_{\sigma}$.
It follows from a certain elementary lemma \cite[V. Lemma 5.1]{hartshorne} that $X \setminus U_{\sigma}$ is of codimension $\ge 2$. 
Therefore, we have: 

\begin{note} 
For any hyperplane $H \in V_P \setminus (\gamma(X_{\rm sm}) \cup S_X \cup \Sigma_P)$, $X_H \cap U_{\sigma}$ is non-empty. 
\end{note} 

We also define $U_G:=\bigcap_{\sigma \in G_P} U_{\sigma}$. 
We have the following natural property in Galois extension (cf. \cite[III. 7.1]{stichtenoth}):

\begin{lemma} 
Let $P$ be a Galois point. 
For any points $Q \in U_G$ and $R \in U_P$ with $\pi_P(Q)=\pi_P(R)$, there exists an element $\sigma \in G_P$ such that $\sigma(Q)=R$. 
\end{lemma} 

\begin{proof}
Assume that $\sigma(Q) \ne R$ for any $\sigma \in G_P$. 
We take a function $f \in \bigcap_{\sigma} \mathcal{O}_{\sigma(Q)} \cap \mathcal{O}_R$ such that $f(\sigma(Q)) \ne 0$ for any $\sigma \in G_P$ and $f(R)=0$. 
Then, we consider the function $g=\prod_{\sigma} \sigma^{*}f$. 
For any $\sigma \in G_P$, we have $\sigma^{*}g=g$. 
Therefore, we find that $g \in \mathcal{O}_{\pi_P(R)}$. 
Note that $g \in m_{\pi_P(R)}$ because of $\pi_P^*m_{\pi(R)}=m_R$ by $R \in U_P$, where $m_R$ is the maximal ideal of the local ring of $\mathcal{O}_R$.  
This implies that $0 \ne g(\pi_P(Q))=g(\pi_P(R))= 0$. 
This is a contradiction. 
\end{proof}

\begin{proof}[Proof of Theorem \ref{HyperplaneSection}] 
Assume that $P$ is Galois. 
We have Theorem \ref{HyperplaneSection} (i) by Lemma 1. 
Let $d_0$ be the degree of the function field extension $K(X)/\pi_P^*K(\Bbb P^n)$. 
Let $H$ be any hyperplane passing through $P$ and satisfying the condition $(\star)$. 
Then, we consider a group morphism
$$\phi(H, P): G_P(X) \rightarrow G; \sigma \mapsto \sigma|_{X_H} $$
where $G=\{ \sigma \in {\rm Bir}(X_H)|\sigma(X_H \cap l) \subset X_H \cap l \mbox{ for a general line } l \mbox{ such that } P \in l \subset H  \}$. 
It follows from Note 5 that $\sigma$ is defined and isomorphic at a general point of $X_H$. 
We also find that $\sigma(X_H)=X_H$. 
Therefore, $\phi(H, P)$ is well-defined.  
Now we prove that $\phi(H,P)$ is an injection. 
We assume that $\sigma|_{X_H}={\rm id}_{X_H}$ and $\sigma \ne {\rm id}_X$. 
Then, $\sigma(Q)=Q$ for any $Q \in X_H$. 
It follows from Lemma 2 that the cardinality $X_H \cap l\setminus{P}$ is at most $d_0-1$ for a general line $l$ in $H$ containing $P$. 
If $\dim X_H=1$, then $l$ is a tangent line of $X_H$ by B\'{e}zout theorem. 
However, this is a contradiction with $(\star)$.  
If $\dim X_H >1$, then this is also a contradiction by taking some linear section of $X_H$ and using Lemma 1. 
Since the order of $G$ is at most $d_0$, $\phi$ is an isomorphism. 
Therefore, $P$ is a Galois point for $X_H$. 
We have (ii) and (iii). 
\end{proof}

\begin{remark} 
It seems to be already known that the assertion in Theorem \ref{HyperplaneSection} holds when $p=0$, $s(X) \le 0$ and $H \in V_P$ is assumed to be general in (ii) (cf. \cite{takahashi}, \cite[Proposition 2.5]{yoshihara2}, \cite{yoshihara3}). 
Yoshihara \cite{yoshihara5} also informed the first author an idea of the proof in this case. 
\end{remark}

\section{Galois points for a normal hypersurface in characteristic zero} 
In this section, we assume that $p=0$.

\subsection{Preliminaries}
Let $X \subset \mathbb{P}^{n+1}$ be an irreducible and reduced hypersurface of degree $d$ and $P$ a point in $\mathbb{P}^{n+1}$.
If $P \not\in X$, then we define the multiplicity of $X$ at $P$ as $0$.  
Let $(X_0: X_1: \dots: X_{n+1})$ be a system of homogeneous coordinates, $(x_1, \dots, x_{n+1})=(X_1/X_0, X_{1}/X_{0}, \dots, X_{n+1}/X_{0})$ a system of affine coordinates, and we assume that $P=(0,0,\dots, 0)$. Let $F(X_0, X_1, \dots, X_{n+1})=0$ be the defining homogeneous equation of $X$, and $f(x_1, \dots, x_{n+1})$ $=F(1, x_1, \dots, x_{n+1})$ $= \sum_{i=0}^{d}f_i$ its dehomogenized polynomial, where $f_i$ is a homogeneous part of $f$ with degree $i$. 
The multiplicity is equal to $m$ if and only if $f_m \ne 0$ and $f_i=0$ for any $i < m$. 

\begin{lemma} \label{l1} 
Let $0 \le m \le d-1$. 
The following three conditions are equivalent{\rm :}
\begin{itemize}
\item[(1)$_m$] The multiplicity at $P$ is $m$, $f_m$ divides $f_{m+1}$, 
and $$f_{m+i}=\binom{d-m}{i} f_m \left(\frac{f_{m+1}}{(d-m)f_m}\right)^i  \,\,\, (i=0,1,\dots, d-m-1)$$
\item[(2)$_m$] By taking a suitable projective transformation fixing the point $P$, the defining equation can be given by 
$$g_m(x_1, \dots, x_{n+1}) + g_d(x_1, \dots, x_{n+1})=0,$$
where $g_m$ and $g_d$ are homogeneous polynomials of degree $m$ and degree $d$, respectively. 
\item[(3)$_m$] The point $P$ is of multiplicity $m$ and Galois, and the birational map induced by $\sigma$ is a restriction of a projective transformation of $\mathbb{P}^{n+1}$ for any $\sigma \in G_P(X)$. 
\end{itemize}
Furthermore, $P$ is Galois and the Galois group $G_P(X)$ is a cyclic group of order $d-m$, if one of the conditions holds. 
\end{lemma} 

\begin{proof}
First let us prove the implication $(1)_m \Rightarrow (2)_m$. 
Since $f_{m+1}/((d-m)f_m)$ is a homogeneous polynomial of degree one, let us put $h:=f_{m+1}/((d-m)f_m)=m_{1}x_1+m_{2}x_2+ \cdots + m_{n+1}x_{n+1}$, where $m_{1}, \cdots, m_{n+1} \in K$. 
Let $\hat{X}_0=X_0+h(X_1, \ldots, X_{n+1})$. 
Then,
$F(\hat{X}_0-h, X_1, \cdots, X_{n+1})=f_m(\hat{X}_0-h)^{d-m}+\cdots+f_{m+i}(\hat{X}_0-h)^{d-m-i}+\cdots+f_d$. 
The coefficient $F_{d-m-i}$ of $\hat{X}_0^{d-m-i}$ is 
$$ \sum_{k=0}^i \binom{d-m-k}{d-m-i} f_{m+k}(-h)^{i-k}. $$
By using $(1)_m$, we have 
$$ F_{d-m-i}=f_mh^i\sum_{k=0}^i(-1)^{i-k}\binom{d-m-k}{d-m-i} \binom{d-m}{k} = (-1)^i f_mh^i \binom{d-m}{d-m-i} \sum_{k=0}^i(-1)^{k} \binom{i}{k}. $$ 
Since $\sum_{k=0}^i(-1)^{k} \binom{i}{k}=0$ if $i>0$, we have the assertion of $(2)_m$. 

Now we prove $(2)_m \Rightarrow (1)_m$. 
Let $\phi$ be the projective transformation as in the condition $(2)_m$ and let $\psi$ be the inverse, and let $\hat{X}_i=\psi^*X_i$ for $0 \le i \le n+1$. 
Then, by the assumption, $g_m(\hat{X}_1, \ldots, \hat{X}_{n+1})\hat{X}_0^{d-m}+g_d(\hat{X}_1, \ldots, \hat{X}_{n+1})=0$. 
Let $A_{\psi}$ be a matrix representing $\psi$. 
Since $\psi(P)=P$, we may assume 
$$ A_{\psi}= \begin{pmatrix} 1 & 0 & \cdots & 0 \\
                   a_{1,0} & a_{1,1} & \cdots & a_{1, n+1} \\
                   \vdots & \vdots & & \vdots \\
                   a_{n+1, 0}& a_{n+1, 1} & \cdots & a_{n+1, n+1} 
\end{pmatrix}.$$
Let $\hat{g}_m(X_0,\dots, X_{n+1}):=g_m(\hat{X}_0, \ldots, \hat{X}_{n+1})$, let $\hat{g}_d(X_0, \ldots, X_{n+1}):=g_d(\hat{X}_0, \ldots, \hat{X}_{n+1})$ and $h:=a_{1, n+1}X_1+\dots+a_{n+1, n+1}X_{n+1}$. 
Then, we have $F=\hat{g}_m(h+X_0)^{d-m}+\hat{g}_d$ and 
$$ f_{m+i}=\hat{g}_m \binom{d-m}{i} h^i.  $$ 
Therefore, $f_{m+1}/f_m=(d-m)h$ and 
$$ \binom{d-m}{i} f_m \left(\frac{f_{m+1}}{(d-m)f_m}\right)^i
= \binom{d-m}{i} f_m h^i = \binom{d-m}{i} \hat{g}_m h^i=f_{m+i} $$
for $i=0, \ldots, d-m-1$. 

The proof of $(2)_m \Leftrightarrow (3)_m$ is the same as the proof for $n=1$ (See \cite[Proposition 1]{yoshihara4}, and also \cite[Proposition 4]{miura}).  
\end{proof}

Below, we assume that $s(X) \le n-2$ and a Galois point is contained in $X_{\rm sm} \cup (\Bbb P^{n+1} \setminus X)$. 
We have the following Proposition, by Lemma \ref{l1} and Theorem \ref{HyperplaneSection}. 

\begin{proposition} \label{l2}
We consider the following weaker condition than $(3)_1$: 
\begin{itemize}
\item[(3')$_1$] The point $P \in X_{\rm sm}$ is Galois. 
\end{itemize}
Then, the three conditions $(1)_1, (2)_1$ and $(3')_1$ are equivalent. 

Especially, if $P \in X_{\rm sm}$ is Galois, then the birational transformation of $X$ induced by $\sigma \in G_P(X)$ is a restriction of a projective transformation of $\Bbb P^{n+1}$. 
\end{proposition}

\begin{proof}
The implications  $(1)_1 \Leftrightarrow (2)_1$ and $(2)_1 \Rightarrow (3')_1$ are clear from Lemma~\ref{l1}. 

Let us prove the implication $(3')_1 \Rightarrow (1)_1$ by induction on the dimension $n$. 
If $n=1$, then it holds by \cite[Proposition~5]{yoshihara1} and the implication of $(2)_1 \Rightarrow (1)_1$ in Lemma~\ref{l1}.  
We assume that $n \geq 2$ and $P$ is a smooth Galois point. Let $H$ be a general hyperplane given by the equation $x_1= a_2 x_2 + \dots + a_{n+1} x_{n+1}$, where $a_i \in K$. 
We put $\tilde{x} = a_2 x_2 + \dots + a_{n+1} x_{n+1}$. 
Then, from Theorem~\ref{HyperplaneSection}, $X_H := X \cap H$ satisfies the condition $(\star)$, and $P$ is a smooth Galois point with respect to $X_P$. 
Hence, by the assumption of the induction, we have that $f_1(\tilde{x}, x_1, \dots, x_{n+1}) \ne 0$, $f_2(\tilde{x}, x_1, \dots, x_{n+1})/f_1(\tilde{x}, x_1, \dots, x_{n+1}) \in K[x_0, \dots, x_{n-1}]$ and 
$$ 
f_{i+1}(\tilde{x}, x_1, \dots, x_{n+1}) 
=\binom{d-1}{i} f_1(\tilde{x}, x_1, \dots, x_{n+1}) \left(\frac{f_2(\tilde{x}, x_1, \dots, x_{n+1})}{(d-1)f_1(\tilde{x}, x_1, \dots, x_{n+1})}\right)^i
$$ 
$(i=0, \dots, d-2)$, for a general $(a_2, \dots, a_{n+1})$. 
This implies the assertion $(1)_1$. 
\end{proof}

\begin{corollary} \label{f_2^2=3f_1f_3}
Assume that $d=4$. Let $P=(1:0:\dots:0)$ be a smooth point of $X$. Then, $P \in \Delta(X)$ if and only if $f_2^2=3f_1f_3$.
\end{corollary}

From Lemma~\ref{l2}, we infer the following remark, which is useful to find smooth Galois points. 
\begin{remark} \label{re1} $ $
\begin{enumerate}
\renewcommand{\labelenumi}{\rm (\arabic{enumi})}
\item If $P \in X_{\rm sm}$ is a Galois point, then $\Bbb T_PX \cap X$ is a (possibly reducible) cone and $P$ is its vertex. 
\item Let $H=H(F)$ be the Hessian of $F$. 
If $P \in X_{\rm sm}$ is a Galois point, then $H(F)(P)=0$. 
\end{enumerate}
\end{remark}

Similar to Proposition \ref{l2}, we have the following: 

\begin{proposition} \label{l2-2}
We consider the following weaker condition than $(3)_0$: 
\begin{itemize}
\item[(3')$_0$] The point $P \in \Bbb P^{n+1} \setminus X$ is Galois. 
\end{itemize}
Then, the three conditions $(1)_0, (2)_0$ and $(3')_0$ are equivalent. 

Especially, if $P \in \Bbb P^{n+1} \setminus X$ is Galois, then the birational transformation of $X$ induced by $\sigma \in G_P(X)$ is a restriction of a projective transformation of $\Bbb P^{n+1}$. 
\end{proposition}

Propositions \ref{l2} and \ref{l2-2} imply that some methods by Yoshihara \cite{yoshihara3} for smooth hypersurfaces are available also for normal hypersurfaces in our situation. 
We introduce the following: 

\begin{definition}
Let $P$ be a Galois point. 
Then, the set $\Bbb F_P:=\{Q \in \Bbb P^{n+1}| \sigma (Q)=Q \mbox{ for any } \sigma \in G_P(X) \} \setminus \{P\}$ is a hyperplane. 
We call $\Bbb F_P$ the fixed hyperplane at $P$. 
\end{definition} 

Note that $\Bbb F_P$ is defined by $X_0=0$ if the defining polynomial has the form in Lemma \ref{l1} (2). 
Now we mention independent Galois points, which is a useful notion to count the number of Galois points (\cite[Definition 4]{yoshihara3}).

\begin{definition}
A set of Galois points $\{P_0, \ldots, P_r\}$ is said to be independent (or, simply, points $P_0, \ldots, P_r$ are said to be independent) if for any two points $P_i$ and $P_j$ ($0 \le i,j \le r$) all the Galois points for $X$ lying on the line $\overline{P_iP_j}$ are exactly $P_i$ and $P_j$. 
\end{definition}

We have the following lemma also for normal hypersurfaces by copying the proof of \cite[Lemma~3]{yoshihara3}. 

\begin{lemma}\label{l5}
If $P_0, \ldots, P_r$ are independent Galois points, then we can choose coordinates $(X_0, \ldots, X_{n+1})$ satisfying $X_j(P_i)=\delta_{ji}$ ($0 \le i \le r$, $0 \le j \le n+1$) and a generator $\sigma_i$ of $G_{P_i}$ ($0 \le i \le r$) has a representation as ${\rm diag}[\zeta, \ldots, \zeta, 1, \zeta, \ldots, \zeta]$, where $1$ is in $i$-th position and $\zeta=e_{d-1}$ (resp. $e_d$). 
Especially we have $r \le n+1$. 
\end{lemma} 

\subsection{Distribution of Galois points for a cone variety}
Firstly, we mention the distribution of Galois points for a cone variety. 
We does not assume the normality of $X$ in this subsection, that is, $s(X)$ may be $n-1$. 
Assume that $X \subset \Bbb P^{n+1}$ is a cone. 
Then, there exists a linear space $M_1$ and $M_2$ in $\Bbb P^{n+1}$ with $M_1 \cap M_2 = \emptyset$ and $M_1 \sharp M_2=\Bbb P^{n+1}$ such that $Y:=X \cap M_1$ is an irreducible hypersurface and $M_2$ is the maximal vertex (i.e. $M_2$ is the vertex with maximal dimension). 
Let $n-a$ be the dimension of $M_2$ (with $0 \le n-a \le n-2$). 
Then, we may assume that $M_1$ is defined by $X_{a+1}=\dots=X_{n+1}=0$, $M_2$ is defined by $X_0=\dots=X_{a}=0$ and $X$ is defined by $F(X_0, \ldots, X_{a})=0$. 

\begin{lemma} \label{Cone1}
Let $P \in M_1$ and $Q \in (P \sharp M_2) \setminus M_2$. 
Then, $P$ is Galois for $X$ if and only if $Q$ is Galois for $X$. 
\end{lemma}

\begin{proof}
By direct computations, $\pi_P$ and $\pi_Q$ induce the same function field extension. 
\end{proof}

\begin{lemma} \label{Cone2}
Let $P \in M_1$. 
Then, $P$ is Galois for $Y$ if and only if $P$ is Galois for $X$.  
\end{lemma}

\begin{proof}
By taking suitable coordinates, we may assume that $P=(1:0:\dots:0)$, and $\pi_P:X \dashrightarrow \Bbb P^n$ is given by $(X_0:X_1:\dots:X_{n+1}) \mapsto (X_1:\dots:X_{n+1})$. Let $(x_0, x_2, \dots, x_{n+1}) = (X_0/X_1, X_2/X_1, \dots, X_{n+1}/X_1)$ be a system of affine coordinates. We have $K(\Bbb P^{a-1})=k(x_2, \dots, x_a)$, $K(\Bbb P^n)=K(\Bbb P^{a-1})(x_{a+1}, \dots, x_{n+1})$, $K(X)=K(\Bbb P^{n})(x_0|_X)$ and $K(Y)=K(\Bbb P^{a-1})(x_0|_Y)$. Let $p_X(x)$ and $p_Y(x)$ be minimal polynomials of $x_0|_X$ over $K(\Bbb P^{n})$ and of $x_0|_Y$ over $K(\Bbb P^{a-1})$, respectively. Since $X$ and $Y$ are given by the same equation $F(X_0, \dots, X_a)=0$, we see that $p_X(x)=p_Y(x)$. Let $L_X$ and $L_Y$ be splitting fields for $p_X(x)$ and $p_Y(x)$, respectively. Then, by \cite[Theorem~29]{artin},  
we have that $\Gal(L_X/K(\Bbb P^{n}))$ is isomorphic to the subgroup of $\Gal(L_Y/K(\Bbb P^{a-1}))$ having $L_Y \cap K(\Bbb P^{n})$ as its fixed field. Here, we note that $L_Y \cap K(\Bbb P^{n}) = K(\Bbb P^{a-1})$, since $x_{a+1}, \dots, x_{n+1}$ are transcendental over $K(\Bbb P^{a-1})$. Namely, we have that $\Gal(L_X/K(\Bbb P^{n}))$ is isomorphic to $\Gal(L_Y/K(\Bbb P^{a-1}))$. Especially, we have $[L_X:K(\Bbb P^{n})]=[L_Y:K(\Bbb P^{a-1})]$. Hence, note that $[K(X):K(\Bbb P^n)]=[K(Y):K(\Bbb P^{a-1})]$, we conclude that $K(X)/K(\Bbb P^n)$ is Galois if and only if $K(Y)/K(\Bbb P^{a-1})$ is Galois. 
\end{proof}

It follows from Lemmas \ref{Cone1} and \ref{Cone2} that we have the following: 
\begin{proposition} \label{Cone3} 
Let $X$ be a cone described as above. 
Then, 
$$ \Delta(X)=(\Delta(Y) \sharp M_2) \setminus M_2 \mbox{ and } \Delta'(X)=(\Delta'(Y) \sharp M_2) \setminus M_2.  $$
In particular, if $X$ is a cone, then $\Delta(X)$ and $\Delta'(X)$ are empty or infinite respectively. 
\end{proposition}

Using this proposition, Theorem \ref{InfiniteGalois} will be proved by Propositions \ref{Infinite1}, \ref{Infinite2} and \ref{Infinite3} which will be proved in the next subsection.

\subsection{Inner Galois points} 

In this subsection, we consider only Galois points contained in $X_{\rm sm}$. 

\begin{lemma} \label{l3}
If $l$ be a line lying on $X$, then the number of Galois points for $X$ on $l$ is zero, one, two or infinitely many. 
The last case occurs only if $X$ is a cone.  
\end{lemma}

\begin{proof} 
Suppose that there exist three Galois points $P_1$, $P_2$ and $P_3$ for $X$ on $l$. 
Then, from Proposition~\ref{l2}, we may assume $P_1=(1:0:\dots:0)$ and 
$X$ is given by the equation 
$$X_1X_{0}^{d-1} + G(X_1, \dots, X_{n+1}) = 0.$$ 
Since $l$ is contained in the tangent space at $P_1$, we may assume that $l$ is given by the equation $X_1=X_3=\dots=X_{n+1}=0$.   
The Gauss map $\gamma$ is given by 
$$((d-1)X_1X_0^{d-2} : X_0^{d-1}+\partial G/\partial X_{1} : \partial G/\partial X_2 : \dots : \partial G/\partial X_{n+1} ).$$
Here, we consider the restriction of $\gamma$ to $l$, which is given by
$$(0: X_0^{d-1}+a_1X_2^{d-1}: da_2X_2^{d-1}:a_3X_2^{d-1}:\dots:a_{n+1}X_{2}^{d-1})$$
where $a_i \in K$ is the coefficient of $X_iX_{2}^{d-1}$ in $G$. 
Assume that the restriction $\gamma|_l$ is not a constant map, i.e., there exists a non-zero coefficient $a_i$ for some $i=0, \dots, n-1$. 
Then, the map $\gamma|_l: l \rightarrow \gamma|_l(l) \cong \mathbb{P}^1$ is  a finite morphism of degree $d-1$. 
It is clear that the Galois point $P_1$ is a ramification point of $\gamma|_l$. So, other Galois point $P_2$ and $P_3$ must be also ramification points. 
However, the number of ramification points of $\gamma|_l$ must be two, this is a contradiction. 
So we have that $\gamma|_l$ is a constant map. Especially, we have that $\Bbb T_{P_1}X=\Bbb T_{P_2}X=\Bbb T_{P_3}X$.

By Remark~\ref{re1}, $\Bbb T_{P_1}X \cap X$ is a cone and $P_1$, $P_2$ and $P_3$ are its vertexes. 
We put that $P_2=(a:0:1:0:\dots:0)$, $P_3=(b:0:1:0:\dots:0)$ ($a,b \in K$) and assume that $a\ne0$. 
Then, calculating the local equation of $\Bbb T_{P_1}X$ at $P_2$, we see that $G(0, 1, X_3, \ldots, X_{n+1})$ is a homogeneous polynomial of degree $d$. 
Hence, we may assume that $X$ is given by the equation 
$$X_1X_{0}^{d-1} + X_2 ( G_0 X_{2}^{d-1} + G_1 X_{2}^{d-2} + \dots + G_{d-2} X_{2}) + G_d = 0,$$ 
where $G_i=G_i(X_0, X_1, X_3, \ldots, X_{n+1})$ $(i=0,1,\cdots, d-2, d)$ is a homogeneous polynomial of degree $i$. 
Examining the condition $(1)_1$ of Lemma~\ref{l1} at $P_2$, we obtain that $G_0 = G_1 = \cdots = G_{d-2} = 0$. Namely, $X$ is given by the equation $$X_1X_{0}^{d-1} + G_d(X_0, X_1, X_3, \ldots, X_{n+1}) = 0,$$and $X$ is a cone with the vertex $O=(0:0:1:0:\dots:0)$. 
\end{proof}

\begin{lemma}\label{l4} Assume that $d \geq 5$. 
If $l$ is a line which does not lie on $X$, then the number of Galois points for $X$ on $l$ is at most one. 
\end{lemma}
\begin{proof} 
Suppose that there exist two Galois points for $X$ on the line $l$. 
Then, we denote them by $P$ and $P'$. 
From Proposition~\ref{l2}, we may assume that $P=(1:0:\dots:0)$ and $X$ is given by the equation $F=X_1X_{0}^{d-1}+G(X_1, \dots, X_{n+1})$. 
Note that if $P' \in X$ is on the hyperplane $X_{0}=0$, then $P \in T_{P'}X$. 
So, we infer from Remark~\ref{re1} that $P'$ is not on the hyperplane $X_{0}=0$. 
Hence, we may assume that $P'=(1:1:0\dots:0)$ and $l$ is given by the equations $X_2=\dots=X_{n+1}=0$. 
Then, the local equation of $X$ at $P'$ is the following: 
$$(w+1)^{d-1}+G(1, u_2, \dots, u_{n+1})=0,$$
where $(w, u_2, \dots, u_{n+1})=(X_{0}/X_1-1,X_2/X_1, \dots, X_{n+1}/X_1)$ is a system of local coordinates. 
Putting $G(1, u_2, \dots, u_{n+1}) = \sum_{i=0}^d g_i$, where $g_i = g_i(u_2, \dots, u_{n+1})$ $(i=0,1,\cdots, d)$ is a homogeneous polynomial of degree $i$. 
Examining Condition $(1)_1$ in Lemma~\ref{l1}, we obtain that 
$$\binom{d-1}{3}{w}^3 +g_3= \binom{d-1}{2} \left( (d-1)w+g_1 \right) \left(\frac{\binom{d-1}{2}{w}^2+g_2}{(d-1)\left( (d-1)w+g_1 \right) }\right)^2$$
Comparing the coefficients of $w^3$ of both sides, then we have a contradiction, if $d \ge 5$. 
\end{proof}

\begin{proposition} \label{Infinite1} 
If $d \ge 5$ and $\Delta(X)$ is infinite, then $X$ is a cone. 
\end{proposition}

\begin{proof}
Suppose that $\Delta(X)$ is an infinite set. 
Then, we infer from Lemma~\ref{l5} that there exist three smooth Galois points $P_1$, $P_2$, $P_3$ which are collinear. 
By Lemma~\ref{l4}, the line $l$ passing through $P_1$, $P_2$, $P_3$ is contained in $X$. So, from Lemma~\ref{l3}, $X$ is a cone. 
\end{proof} 

On the other hand, if $d \ge 5$, and if $\Delta(X)$ is finite and non-empty, then the set $\Delta(X)$ is independent by Proposition \ref{Cone3} and Lemmas \ref{l3} and \ref{l4}. 
Theorem 1 (II) is derived from the following lemma which is a generalization of \cite[Lemma~4]{yoshihara3}: 
 
\begin{lemma}  \label{IndependentInner}
Let $m=[(n+s+1)/2]$. 
The cardinality of a set of independent inner Galois points is at most $m+1$. 
The equality holds if and only if $X$ is projectively equivalent to the hypersurface defined by 
$$X_{m+1}X_0^{d-1}+\dots+X_{2m-s}X_{m-s-1}^{d-1}+A_{m-s}X_{m-s}^{d-1}+\dots+A_mX_m^{d-1}+G=0 $$
where $A_{m-s}, \ldots, A_m , G \in K[X_{m+1}, \ldots, X_{n+1}]$ are homogeneous polynomials with $\deg A_{m-s}=\dots=\deg A_m=1$ and $\deg G=d$.  
(Note that the polynomial $A_{m-s}X_{m-s}^{d-1}+\dots+A_{m}X_m^{d-1}$ does not appear if $s=-1$.)  
\end{lemma}

\begin{proof}
Let $P_0, \ldots, P_r$ form a set of independent inner Galois points. 
Suppose that $r \ge m+1$. 
Then take a system of coordinates $(X_0, \ldots, X_{n+1})$ satisfying that $X_j(P_i)=\delta_{ji}$, where $0 \le i \le m+1$ and $0 \le j \le n+1$. 
By Lemma \ref{l5}, we can assume that $\sigma_i$ is a diagonal matrix ${\rm diag}[\zeta, \ldots, \zeta, 1, \zeta, \ldots, \zeta]$, where $\zeta=e_{d-1}$. 
Since $F^{\sigma}=\lambda_iF$ for $\lambda_i \in K \setminus 0$, we infer that $F$ has the expression as 
$$ F=A_0X_{0}^{d-1}+\dots+A_{m+1}X_{m+1}^{d-1}+G, $$
where $A_i$ and $G$ are forms in $K[X_{m+2}, \ldots, X_{n+1}]$, and $\deg A_i=1$ and $\deg G=d$. 
Let $t:=\dim \langle A_0, \ldots, A_{m+1} \rangle-1$. 
Then, $0 \le t \le n-m-1 \le m$. 
We may assume that $A_0, \ldots, A_{t}$ be a basis of the above vector space and $A_0=X_{m+2}, \ldots, A_{t}=X_{m+t+2}$. 
Then, $A_{t+1}$, \ldots, $A_{m+1}$ are represented by linear combinations of $X_{m+2}, \ldots, X_{m+t+2}$.   
We consider the locus $\Gamma$ which is defined by $X_{m+2}=\dots=X_{m+t+2}=\frac{\partial F}{\partial X_{m+2}}=\dots=\frac{\partial F}{\partial X_{n+1}}=0$. 
We find easily that the locus $\Gamma$ is contained in the singular locus of $X$. 
The dimension of $\Gamma$ is at least $(n+1)-(t+1+n-m)=m-t$. 
Now, $m-t \ge 2m-n+1 \ge s+1$. 
This is a contradiction. 

Now assume that $r=m$. 
Similar to the above discussion, 
$F$ has the expression as 
$$ F=A_0X_{0}^{d-1}+\dots+A_mX_{m}^{d-1}+G, $$
where $A_i$ and $G$ are forms in $K[X_{m+1}, \ldots, X_{n+1}]$, and $\deg A_i=1$ and $\deg G=d$. 
Let $t=\dim \langle A_0, \ldots, A_m \rangle-1$. 
Then, $0 \le t \le \max\{n-m, m\}$. 
We may assume that $A_0, \ldots, A_{t}$ be a basis of the above vector space and $A_0=X_{m+1}, \ldots, A_{t}=X_{m+t+1}$. 
We consider the locus $\Gamma$ which is defined by $X_{m+1}=\dots=X_{m+t+1}=\frac{\partial F}{\partial X_{m+1}}=\dots=\frac{\partial F}{\partial X_{n+1}}=0$. 
We find easily that the locus $\Gamma$ is contained in the singular locus of $X$. 
The dimension of $\Gamma$ is at least $(n+1)-(t+1+n-m+1)=m-t-1$. 
Therefore, $s \ge m-t-1$. 
We have $t \ge m-s-1$. 
Then, we have the assertion. 
\end{proof} 

Now, we consider the case where $d=4$. 

\begin{lemma}\label{OutsideTangent} 
Assume that $d=4$. Let $P \in X$ be a Galois point, $\sigma$ a generator of $G_P$. 
\begin{enumerate}
\item  If a smooth point $Q$ on $X \setminus \mathbb{T}_PX$ is Galois then $\sigma(Q) \ne Q$. 
\item Let $P_1:=P, P_2, P_3, P_4$ be distinct points in $\Delta(X)$. 
Further, we assume that these points are collinear and the line $l$ passing through these points is not contained in $X$. 
Then, by taking a suitable projective transformation, we can express $X$ as $X_1X_{0}^3+X_1^4+H(X_2,\dots,X_{n+1})=0$, where $H(X_2,\dots,X_{n+1})$ is a homogeneous polynomial of degree four. 
\item Assume the same assumptions in (2) and that $X$ is not a cone. 
Then, we have that $\Delta(X) \setminus \mathbb{T}_{P}X = \{P_2, P_3, P_4 \}$. 
\end{enumerate}
\end{lemma}
\begin{proof}
Let us prove the assertion (1). 
Suppose that $\sigma(Q)=Q$. 
By the assumtion, the line $\overline{PQ} \not \subset X$. 
Therefore, $P$ and $Q$ are independent Galois points. 
By Proposition~\ref{l2}, we may assume that $P=(1:0: \dots :0)$ and $Q=(0:1:0:\dots:0)$, and $X$ is given by $A_0X_{0}^3+A_1X_1^3+G(X_1,\dots,X_{n+1})=0$ where $A_0, A_1 \in K[X_2, \ldots, X_{n+1}]$. 
The tangent space at $P$ is defined by $A_0=0$ and $A_0(Q)=0$. 
This contradicts the assumption. 

Let us prove the assertion (2). 
By Proposition~\ref{l2}, we may assume that $P_1=(1:0:\dots:0)$, $X$ is given by the equation $X_1X_{0}^3+G(X_1,\dots,X_{n+1})=0$. 
The point $P_2$ is not contained in $\Bbb T_{P_1}X$ nor $\Bbb F_P$. 
Therefore, we may assume that $P_2=(a:1:0:\dots:0)$ and $l$ is given by the equations $X_2= \dots = X_{n+1}=0$. 
Let $G(X_1,\dots,X_{n+1})=G_0 X_1^4 + G_1 X_1^3 + G_2 X_1^2 + G_3 X_1+ G_4$
where $G_i=G_i(X_2,\dots,X_{n+1})$ is a homogeneous polynomial of degree $i$ and $G_0=-a^3$. 
Let $w:=X_0/X_1-a$, $u_2:=X_2/X_1$, \ldots, $u_{n+1}:=X_{n+1}/X_1$. 
Then, $$F(w+a,1,u_2, \ldots, u_{n+1})=(w+a)^3+(-a)^3+G_1+G_2+G_3+G_4=0. $$
By applying Corollary~\ref{f_2^2=3f_1f_3} at $P_2$, we have 
$$ (3a^2+G_2)^2=3(3a^2w+G_1)(w^3+G_3) $$
as a polynomial. 
This implies that $G_1=G_2=G_3=0$ and hence, the assertion.

Let us prove the assertion (3). 
By the assertion in (2), we may assume that $P_1=(1:0: \dots :0)$, $X$ is given by $X_1X_{0}^3+X_1^4+H(X_2,\dots,X_{n+1})=0$, and $l$ is given by the equations $X_2= \dots = X_{n+1}=0$. 
Suppose that there exists another point $Q$ in $\Delta(X) \setminus \mathbb{T}_{P_1}X$. 
By taking a suitable projective transformation, we may assume that $Q=(a:1:1:0:\dots:0)$. We put that 
$$H(X_2,\dots,X_{n+1})=H_0 X_{2}^4 + H_1 X_{2}^3 + H_2 X_{2}^2 + H_3 X_{2}+ H_4,$$
where $H_i=H_i(X_3, \dots,X_{n+1})$ is a homogeneous polynomial of degree $i$, and $H_0=-a^3-1$. 
Let $w_0:=X_0/X_2-a$, $w_1:=X_1/X_2-1$, $u_3=X_3/X_2, \ldots, u_{n+1}:=X_{n+1}/X_2$. 
Then, $$ F(w_0+a, w_1+1, 1, u_3, \ldots, u_{n+1})=(w_1+1)(w_0+a)^3+(w_1+1)^4+(-a^3-1)+H_1+H_2+H_3+H_4=0. $$
By applying Corollary~\ref{f_2^2=3f_1f_3} at $P_2$, we have 
$$ (a^3w_1+3a^2w_0+4w_1+H_1)(3aw_0^2w_1+w_0^3+4w_1^3+H_3)=3(3a^2w_0w_1+3aw_0^2+6w_1^2+H_2)^2 $$
as a polynomial. 
This implies that $H_1=H_2=H_3=0$ and $a^3=-1$. 
Therefore, $F=X_1X_0^3+X_1^4+H_4$ with $H_4 \in K[X_3, \ldots, X_{n+1}]$ and $X$ is a cone. 
\end{proof}

\begin{proposition} \label{Infinite2} 
If $d=4$ and $\Delta(X)$ is infinite, then $X$ is a cone. 
\end{proposition}

\begin{proof} 
Assume that $X$ is not a cone. 
It follows from Lemma \ref{l3} that if a line $l$ contains three Galois points, then $l \not\subset X$.  
Let $P \in \Delta(X)$.  
Then, we have $\Delta(X) \setminus \Bbb F_P \subset \Delta(X) \setminus \Bbb T_{P}X$, because the line $\overline{PQ}$ contains four Galois points if $Q \in \Delta(X) \setminus \Bbb F_P$. 
It follows from Lemma \ref{OutsideTangent} (3) that $\Delta(X) \setminus \Bbb F_P$ is a finite set. 

Let $P_1$ be a point in $\Delta(X)$. 
Since $\Delta(X) \setminus \Bbb F_{P_1}$ is finite, we can take a point $P_2$ in $\Delta(X) \cap \Bbb F_{P_1}$. 
Then, we infer that the set $\Delta(X) \setminus (\Bbb F_{P_1} \cap \Bbb F_{P_2})=(\Delta(X) \setminus \Bbb F_{P_1}) \cup (\Delta(X) \setminus \Bbb F_{P_2})$ is finite and the set $\Delta(X) \cap (\Bbb F_{P_1} \cap \Bbb F_{P_2} )$ is infinite. 
Hence, we can take a point $P_3$ in $\Delta(X) \cap (\Bbb F_{P_1} \cap \Bbb F_{P_2} )$. 
In this way, we can take infinitely many points $P_1, P_2, \dots$ such that $P_i \in \Delta(X) \cap (\Bbb F_{P_1} \cap \dots \cap \Bbb F_{P_{i-1}})$. 
Note that $\dim(\Bbb F_{P_1} \cap \dots \cap \Bbb F_{P_i}) = \dim(\Bbb F_{P_1} \cap \dots \cap \Bbb F_{P_{i-1}})-1$, because $P_i \not\in \Bbb F_{P_i}$. 
We have that $\Bbb F_{P_1} \cap \dots \cap \Bbb F_{P_{n+2}} = \emptyset$, thus we have a contradiction. 
\end{proof}

Now we consider the case where $\Delta(X)$ is finite and non-empty. 
Then, $X$ is not a cone by Proposition \ref{Cone3} in Subsection 3.2. 
Below in this subsection, 
we define $$r=r(X):=\max\sharp \{P_0, \ldots, P_b: \mbox{ independent Galois points on $X$}\}-1$$ 
and let $P_0, \ldots, P_r$ form a set of independent Galois points. 
(Then, $r \le m$ and $r <m$ is possible.)
We also define 
$$\mu:=\sharp\{l: \mbox{line }| l \mbox{ contains four Galois points } \}-1.$$
We have that $-1 \le \mu \le r$ by the following Lemma:

\begin{lemma}
Let $P_0, \ldots, P_r$ be independent Galois points and let $Q \in X$ be a Galois point distinct from $P_0, \ldots, P_r$. 
Assume that $X$ is not a cone. 
Then,
\begin{itemize} \label{LinesFourGalois}
\item[(1)] There exists $i$ such that $\sharp (\overline{P_iQ} \cap \Delta(X))=4$.
\item[(2)] If $l_1, l_2$ are lines such that $Q \in l_i$ and $\sharp (l_i \cap \Delta(X))=4$ for $i=1,2$, then $l_1=l_2$.  
\end{itemize} 
\end{lemma}

\begin{proof}
Assume that the line $\overline{P_iQ}$ consists exactly two Galois points for any $i$. 
Then, $P_0, \ldots, P_r, Q$ are independent Galois points. 
This contradicts the assumption on $r$. 
We have (1). 
The assertion in (2) is derived from Lemma \ref{OutsideTangent} (3). 
\end{proof}

Let $P_0,\ldots, P_{\mu}$ be Galois points with lines $l_0, \ldots, l_{\mu}$, where $P_i \in l_i$ and $l_i$ contains four Galois points. 
Using Corollary \ref{f_2^2=3f_1f_3}, Lemmas \ref{OutsideTangent} and \ref{LinesFourGalois}, for a suitable coordinate, we find that $X$ is defined by  
$$ X_{\mu+1}X_0^3+\dots+X_{2\mu+1}X_{\mu}^3+X_{\mu+1}^4+\dots+X_{2\mu+1}^4+G=0 $$
where $G \in K[X_{2\mu+2}, \ldots, X_{n+1}]$ (cf. \cite[p. 532]{yoshihara3}). 
Then, by direct computations, the singular locus is contained in a linear space defined by $X_0=\cdots=X_{2\mu+1}=0$. 
This implies that $s \ge (n+1)-(2\mu+2)$ and hence $\mu \le (n+s+1)/2-s-1$. 
Therefore, we have the following:

\begin{lemma} \label{restrict-mu}
$\mu \le m-s-1$.  
\end{lemma}

\begin{proof}[Proof of Theorem 1 (I-0)--(I-4)] 
(0) 
Now, the number of inner Galois points is $4(\mu+1)+(r-\mu)$. 
Note that $\mu \le r \le m$ and $\mu \le m-s-1$, by Lemmas \ref{IndependentInner} and \ref{restrict-mu}. 
Therefore, we have an upper bound $4(m-s)+(s+1)$.

(1) We consider the case where $n+s$ is odd. 
Then, $2m-s=n+1$. 

(1-1) Assume that $\mu=m-s-1$.  
Then, $X$ is projectively equivalent to the hypersurface defined by 
$$X_{m-s}X_0^3+\dots+X_{2m-2s-1}X_{m-s-1}^3+X_{m-s}^4+\dots+X_{2m-2s-1}^4+G=0$$ where $G \in K[X_{2m-2s}, \ldots, X_{n+1}]$. 
By direct computations, the singular locus of $X$ is defined by
$$ X_0=\dots=X_{2m-2s-1}=\frac{\partial G}{\partial X_{2m-2s}}=\dots=\frac{\partial G}{\partial X_{n+1}}=0. $$
We have 
$$ \frac{\partial G}{\partial X_{2m-2s}}=\dots=\frac{\partial G}{\partial X_{n+1}}=0 $$
as a polynomial by counting the dimension of the singularity (and that $X$ is not a cone). 
Then, $G=0$ as a polynomial. 
This implies $s=-1$. 

(1-2) Assume that $s \ge 0$, $\mu=m-s-2$ and $r=m$. 
Then, $X$ is projectively equivalent to the hypersurface defined by \\
 $X_{m+1}X_0^3+\dots+X_{2m-s-1}X_{m-s-2}^3+A_{m-s-1}X_{m-s-1}^3+\dots+A_mX_m^3$ \\
\hspace{80mm} $+X_{m+1}^4+\dots+X_{2m-s-1}^4+G =0 $ \\
where $A_{m-s-1}, \ldots, A_m, G \in K[X_{2m-s}]$. 
Therefore, we may assume that $A_{m-s-1}=\dots=A_m=X_{2m-s}$ and $G=aX^4$ where $a=0$ or $1$. 
Then, $\sharp\Delta(X)=4(\mu+1)+(m-\mu)=4(m-s-1)+(s+2)=4(m-s)+(s-2)$. 

(2) We consider the case where $n+s$ is even.

(2-1) Assume that $\mu=m-s-1$. 
Then $X$ is projectively equivalent to the hypersurface defined by 
$$X_{m-s}X_0^3+\dots+X_{2m-2s-1}X_{m-s-1}^3+X_{m-s}^4+\dots+X_{2m-2s-1}^4+G=0$$  
where $G \in K[X_{2m-2s}, \ldots, X_{n+1}]$. 
Assume that $r \ge \mu+1$. 
Then, $s \ge 0$ and we can take $G=X_{2m-2s+1}X_{2m-2s}^3+G_0$ where $G_0 \in K[X_{2m-2s+1}, \ldots, X_{n+1}]$. 
The singular locus is defined by 
$$ X_{0}=\dots=X_{2m-2s-1}=X_{2m-2s}X_{2m-2s+1}$$ 
$$=X_{2m-2s}^3+\frac{\partial G_0}{\partial X_{2m-2s+1}}=\frac{\partial G_0}{\partial X_{2m-2s+2}}=\dots=\frac{\partial G_0}{\partial X_{2m-2s+1}}=0. $$
Then, by counting the dimension of the singularity of $X$, we should have $G_0=0$ as a polynomial and hence, $s=0$ and $r=\mu+1=m$. 
Then, $\sharp\Delta(X)=4m+1$. 
On the other hand, if $s \ne 0$, then $r = \mu$  and $\sharp\Delta(X)=4(m-s)$. 
Then, by counting the dimension of the singularity of $X$, $G$ should have a multiple component. 

Using the result in (1-1), we find that the bound $4(m-s)+(s+1)$ is sharp if and only if $s=-1$, or $s=0$ and $n$ is even. 

(2-2) Assume that $s \ge 1$, $\mu=m-s-2$ and $r=m$. 
Then, $X$ is projectively equivalent to the hypersurface defined by \\
 $X_{m+1}X_0^3+\dots+X_{2m-s-1}X_{m-s-2}^3+A_{m-s-1}X_{m-s-1}^3+\dots+A_mX_m^3$ \\
\hspace{80mm} $+X_{m+1}^4+\dots+X_{2m-s-1}^4+G =0 $ \\
where $A_{m-s-1}, \ldots, A_m, G \in K[X_{2m-s}, X_{n+1}]$. 
Then, $\sharp\Delta(X)=4(\mu+1)+(m-\mu)=4(m-s-1)+(s+2)=4(m-s)+(s-2)$. 
There exist examples of hypersurfaces such that $s \ge 1$, $\mu=m-s-2$ and $r=m$, as in Example \ref{n+s:even} below. 
\end{proof}

We define $t:=n-\dim \bigcap_{0 \le i \le r}\Bbb T_{P_i}X$. 
We consider examples when $n+s$ is even. 

\begin{example} \label{n+s:even}
Let $s \ge 1$ and $n+s$ is even. 
\begin{itemize}
\item[(i)] A hypersurface in $\Bbb P^{n+1}$ defined by \\
$X_{m+1}X_0^3+\dots+X_{2m-s}X_{m-s-1}^3+X_{2m-s+1}X_{m-s}^3+X_{2m-s}(X_{m-s+1}^3+\dots+X_m^3)$  \\
\hspace{90mm} $+X_{m+1}^4+\dots+X_{2m-s-1}^4=0$ \\
satisfies that $\dim {\rm Sing}(X)=s$, $r=m$, $t=m-s$ and $\mu=m-s-2$. 
\item[(ii)] A hypersurface in $\Bbb P^{n+1}$ defined by \\
$X_{m+1}X_0^3+\dots+X_{2m-s}X_{m-s-1}^3+X_{2m-s}(X_{m-s}^3+\dots+X_m^3)$ \\
\hspace{70mm} $+X_{m+1}^4+\dots+X_{2m-s-1}^4+X_{2m-s+1}^4=0 $ \\
satisfies that $\dim {\rm Sing}(X)=s$, $r=m$, $t=m-s-1$ and $\mu=m-s-2$. 
\end{itemize}
\end{example}

\subsection{Outer Galois points}

\begin{lemma} \label{OuterCollinearGaloisPoints}
Let $l \subset \Bbb P^{n+1}$ be a line. 
Then, the cardinality of the set $\Delta'(X) \cap l$ is zero, one, two or infinity. 
The last case occurs only if $X$ is a cone. 
\end{lemma} 

\begin{proof}
Let $P_1, P_2, P_3 \in \Delta'(X) \cap l$ be distinct three outer Galois points and let $\sigma \in G_{P_1}$ be a generator. 
From Proposition~\ref{l2-2}, we may assume that $P_1=(1:0:\dots:0)$, $X$ is defined by the equation 
$$X_{0}^d+G(X_1,\dots,X_{n+1})=0,$$
where $G(X_1,\dots,X_{n+1})$ is a homogeneous polynomial of degree $d$, 
and $\sigma=\mathrm{diag}[e_d, \dots, e_d, 1]$. Further, we may assume that $l$ is defined by $X_2 = \dots =X_{n+1}=0$. 
Note that $l \cap \Bbb F_{P_1}$ consists only one point. 
Therefore, $P_2$ or $P_3$ is not contained in $\Bbb F_{P_1}$ at $P_1$. 
We may assume that $P_2$ is so and $P_2=(1:1:0:\dots:0)$, since $\Bbb F_{P_1}$ is given by $X_{0}=0$. 
Then, the local equation of $X$ at $P_2$ is the following: 
$$(w+1)^d+G(1, u_2,\dots,u_{n+1}),$$
where $(w, u_2,\dots,u_{n+1})=(X_0/X_1-1, X_2/X_1, \dots, X_{n+1}/X_1)$ is a system of local coordinates. 
We put that $G(1, u_2,\dots,u_{n+1})=\sum_{i=0}^d g_i$, where $g_i=g_i(u_2,\dots,u_{n+1})$ is a homogeneous polynomial of degree $i$. 
Examining the condition $(1)_0$ in Lemma~\ref{l1}, we obtain that 
$$g_i= \binom{d}{i} \left( \frac{(dw+g_1)^i}{d^i(1+g_0)^{i-1}} - w^i \right), \, i=1,\dots, d-1. $$
Calculating the coefficients of $w$ and $w^2$ in the equation where $i=2$, we see that $g_0=0$ and $g_1=0$. Hence, we see that $g_0=\dots=g_{d-1}=0$. 
Thus, we conclude that $X$ is defined by the following: 
$$X_{0}^d+g_d(X_2,\dots,X_{n+1})=0,$$
and $X$ is a cone with a vertex $Q=(0:1:0:\cdots:0)$. 
\end{proof} 

\begin{proposition} \label{Infinite3} 
If $\Delta'(X)$ is infinite, then $X$ is a cone. 
\end{proposition}

\begin{proof}
Assume that $\Delta'(X)$ is infinite. Then, by lemma~\ref{l5}, there exists three outer Galois points which are collinear. By Lemma~\ref{OuterCollinearGaloisPoints}, we have that $X$ is a cone. 
\end{proof}

\begin{lemma}  \label{IndependentOuter} 
The cardinality of a set of independent outer Galois points is at most $n-s+1$.  
The equality holds if and only if $X$ is projectively equivalent to the hypersurface defined by 
$$ X_0^d+\dots+X_{n-s}^d=0.  $$
Therefore, $X$ is smooth or a cone. 

Furthermore, if $s \ge 0$ and $X$ is not a cone, then the cardinality of a set of independent outer Galois points is at most $n-s$. 
Then, the equality holds if and only if $X$ is projective equivalent to the hypersurface defined by 
$$ X_0^d+\dots+X_{n-s-1}^d+G=0 $$
where $G \in K[X_{n-s}, \ldots, X_{n+1}]$ has a multiple component. 
\end{lemma}

\begin{proof} 
Let $P_0, \ldots, P_r$ are independent outer Galois points. 
Then take a system of coordinates $(X_0, \ldots, X_{n+1})$ satisfying that $X_j(P_i)=\delta_{ji}$, where $0 \le i \le m+1$ and $0 \le j \le n+1$. 
By Lemma \ref{l5}, we can assume that $\sigma_i$ is a diagonal matrix ${\rm diag}[\zeta, \ldots, \zeta, 1, \zeta, \ldots, \zeta]$, where $\zeta=e_{d-1}$. 
Since $F^{\sigma}=\lambda_iF$ for $\lambda_i \in K \setminus 0$, we infer that $F$ has the expression as 
$$ F=X_{0}^{d}+\dots+X_{r}^{d}+G, $$
where $G \in K[X_{r+1}, \ldots, X_{n+1}]$. 
Then, the singular locus is defined by 
$$ X_0=\dots=X_r=\frac{\partial G}{\partial X_{r+1}}=\dots=\frac{\partial G}{\partial X_{n+1}}=0. $$
Therefore, $s \le n+1-(r+1)$. 
The equality holds if and only if $G=0$ as a polynomial, hence $X$ is smooth or a cone. 

Assume that $s \ge 0$ and $X$ is not a cone. 
Then, $r \le n$, there exists $j$ with $r+1 \le j \le n+1$ such that $\frac{\partial G}{\partial X_{j}} \ne 0$ as a polynomial, and hence $s \le n+1-(r+2)$. 
The equality holds if and only if $G$ has a multiple component. 
\end{proof}

\begin{proof}[Proof of Theorem \ref{OuterGalois}]
Assume that $\Delta'(X)$ is finite and non-empty. 
Then, $X$ is not a cone by Proposition \ref{Cone3} in Subsection 3.2, and the set $\Delta'(X)$ is independent by Lemma \ref{OuterCollinearGaloisPoints}. 
Therefore, we have the conclusion by Lemma \ref{IndependentOuter}. 
\end{proof}

\section{Galois points for a Fermat hypersurface of degree $p^e+1$ in $p>0$} 
Let $p>0$ and $q \ge 3$ be a power of $p$. 
We consider the Fermat hypersurface $F_n(q+1) \subset \Bbb P^{n+1}$ of degree $q+1$: 
\begin{equation}
X_0^{q+1}+X_1^{q+1}+\cdots+X_{n+1}^{q+1}=0. 
\end{equation} 
Firstly, we prove the ``if'' part of Theorem \ref{Fermat}. 
Homma gave an elegant proof if $n=1$ by using automorphisms of Hermitian curve (\cite[Claims 1 and 2]{homma}).  
We give an elementary proof.  

\begin{proposition}
If $P \in \Bbb P^{n+1}$ is $\Bbb F_{q^2}$-rational, then $P$ is Galois for $F_n(q+1)$. 
\end{proposition} 

\begin{proof}
Let $P=(1:a_1:\cdots:a_{n+1})$ and $\pi_P=(X_1-a_1X_0:\cdots:X_{n+1}-a_{n+1}X_0)$ be the projection from $P$. 
Then, we have a field extension $K(x_0, \ldots, x_n)/K(x_1, \ldots, x_n)$ with $f(x_0, \ldots, x_n)=(1+a_1^{q+1}+\cdots+a_{n}^{q+1}+a_{n+1}^{q+1})x_0^{q+1}+(a_1^qx_1+\cdots+a_{n}^qx_n+a_{n+1}^q)x_0^q+(a_1x_1^q+\cdots+a_nx_n^q+a_{n+1})x_0+(x_1^{q+1}+\cdots+x_n^{q+1}+1)=0$. 

If $P \in X$ and $P$ is $\Bbb F_{q^2}$-rational, then $f(x_0, \ldots, x_n)=(a_1^qx_1+\cdots+a_{n}^qx_n+a_{n+1}^q)x_0^q+(a_1^qx_1+\cdots+a_n^qx_n+a_{n+1}^q)^qx_0+(x_1^{q+1}+\cdots+x_n^{q+1}+1)=0$. 
It is not difficult to check that the extension is Galois. 

We assume that $P \in \Bbb P^{n+1} \setminus X$ and $P$ is $\Bbb F_{q^2}$-rational.  
Let $\beta(x_1, \ldots, x_n)=(a_1^qx_1+\cdots+a_n^qx_n+a_{n+1}^q)/(1+a_1^{q+1}+\cdots+a_{n}^{q+1}+a_{n+1}^{q+1})$ and $\hat{x_0}=x_0+\beta$. 
Then, we have $0=f(x_0, x_1, \ldots, x_n)=f(\hat{x_0}-\beta, x_1, \ldots, x_n)=(1+a_1^{q+1}+\cdots+a_{n}^{q+1}+a_{n+1}^{q+1})\hat{x_0}^{q+1}+g(x_1, \ldots, x_n)$. 
Therefore, our extension is cyclic.   
\end{proof} 

Next, we prove the ``only-if'' part. 
We will use the result of Homma \cite[Claim 3]{homma} with $n=1$, because we will use induction on the dimension $n$.

\begin{proposition}
If $P \in \Bbb P^{n+1}$ is a Galois point for $F_n(q+1)$, then $P$ is $\Bbb F_{q^2}$-rational. 
\end{proposition}

\begin{proof}
Let $P=(1:a_1:\cdots:a_{n+1})$. 
We may assume that $i_0$ is an integer such that $a_i \ne 0$ and $1+a_i^{q+1} \ne 0$ if and only if $i \le i_0$. 
Note that $a_i \in \Bbb F_{q^2}$ for any $i > i_0$ because $a_i=0$ or $1+a_i^{q+1}=0$. 
Let $1 \le i \le i_0$ and let $H_i$ be a hyperplane defined by $X_i-a_iX_0$. 
Then, $X \cap H_i$ is defined by 
$$ (1+a_i^{q+1})X_0^{q+1}+\cdots+X_{i-1}^{q+1}+X_{i+1}^{q+1}+\cdots+X_{n+1}^{q+1}=0. $$
We take a linear transformation $\phi_i$ on $H_i \cong \Bbb P^{n}$ defined as follows: 
$(1:x_1:\cdots:x_{i-1}:x_{i+1}:\cdots:x_{n+1})$ 
$\mapsto (1:x_1/\sqrt[q+1]{1+a_i^{q+1}}:\cdots:x_{i-1}/\sqrt[q+1]{1+a_i^{q+1}}:x_{i+1}/\sqrt[q+1]{1+a_i^{q+1}}:\cdots:x_{n+1}/\sqrt[q+1]{1+a_i^{q+1}})$.  
Then $\phi_i(X \cap H_i)$ is the Fermat hypersurface in $H_i \cong \Bbb P^{n}$. 
Note that $H_i$ satisfies the condition $(\star)$ in Theorem \ref{HyperplaneSection}. 

We will prove the assertion by induction on the dimension $n$. 
If $n=1$, then the assertion holds by a result of Homma \cite{homma}. 
We assume that $n \ge 2$ and $P$ is a Galois point. 
It follows from Theorem \ref{HyperplaneSection} (ii) that $P$ is also Galois for $X \cap H_i$. 
By the assumption of the induction, $\phi_i(P)$ is a $\Bbb F_{q^2}$-rational point.  
Now we have 
$$ \phi_i(P)=(1:a_1/\sqrt[q+1]{1+a_i^{q+1}}:\cdots:a_{n+1}/\sqrt[q+1]{1+a_i^{q+1}}).$$ 
Therefore, $(a_j/\sqrt[q+1]{1+a_i^{q+1}})^{q^2-1}=1$ for any $i \ne j$ with $1 \le i, j \le i_0$. 
By direct computations, we have $a_i \in \Bbb F_{q^2}$ for any $i$. 
\end{proof}

\section{Examples} 

We give examples of hypersurfaces which have infinitely many Galois points and are not cones, in $p>0$. 

\begin{example} \label{example1} 
Let $X \subset \Bbb P^3$ be a hypersurface defined by 
$$ F=ZW^p-X^pW-Y^{p+1}=0. $$ 
Let $P=(0:0:1:0)$, $Q=(1:0:0:0)$ and let $L$ be a line defined by $Y=W=0$. 
Then, we have the followings: 
\begin{itemize}
\item[(i)] ${\rm Sing}(X)=\{P\}$. 
\item[(ii)] $X$ is not a cone. 
\item[(iii)] $L \setminus \{P, Q\} \subset \Delta(X)$. 
\end{itemize}
\end{example} 

\begin{proof}
Note that $X$ is the closure of the image of a morphism 
$$ \phi:\Bbb A^2 \rightarrow \Bbb P^3; \ (x,y) \mapsto (x:y:x^p+y^{p+1}:1). $$ 
The Gauss map $\gamma$ is given by 
$$ (\partial F/\partial X:\partial F/\partial Y:\partial F/\partial Z:\partial F/\partial W)= (0:-Y^p:W^p:-X^p). $$
We have (i) by the form of the Gauss map. 
We also have (ii) because $\gamma$ is generically finite by that $\gamma \circ \phi=(0:-y^p:1:-x^p)$. 
Now, we prove (iii).  
Let $R=(1:0:a:0)$ and $\pi_R=(Z-aX:Y:W)$ be the projection. 
Then, we have the field extension $K(x,y,z)/K(y,z)$ with $$x^p-ax-z+y^{p+1}=0. $$ 
Therefore, this extension is Galois if $a \ne 0$. 
\end{proof} 

\begin{example}
Let $X \subset \Bbb P^3$ be a hypersurface defined by 
$$ F=ZW^{p^2-1}-X^pW^{p^2-p}-Y^{p^2}=0. $$ 
Let $L_1$ be a line $Y=W=0$, $L_2$ be a line $Z=W=0$ and $H$ be a plane defined by $W=0$. 
Then, we have the followings: 
\begin{itemize}
\item[(i)] ${\rm Sing}(X)=L_1$. 
\item[(ii)] $X$ is not a cone. 
\item[(iii)] $H \setminus (L_1 \cup L_2) \subset \Delta'(X)$. 
\end{itemize}
\end{example} 

\begin{proof}
Note that $X$ is the closure of the image of a morphism 
$$ \phi:\Bbb A^2 \rightarrow \Bbb P^3; \ (x,y) \mapsto (x:y:x^p+y^{p^2}:1). $$ 
The Gauss map $\gamma$ is given by 
$$ (\partial F/\partial X:\partial F/\partial Y:\partial F/\partial Z:\partial F/\partial W)= (0:0:W^{p^2-1}:-ZW^{p^2-2}). $$
We have (i) by the form of the Gauss map. 
Now, we have $$\gamma \circ \phi(x,y)=(0:0:1:-x^p-y^{p^2}). $$  
Therefore, the general fiber of the Gauss map is an irreducible curve defined by $x+y^p=c$ (for a suitable constant $c$). 
Then, we have (ii). 
We prove (iii). 
Let $R=(1:a:b:0)$ and $\pi_R=(Z-aX:Y-bX:W)$ be the projection. 
Then, we have the field extension $K(x,y,z)/K(y,z)$ with $$b^{p^2}x^{p^2}+x^p-ax+z-y^{p+1}=0. $$ 
Therefore, this extension is Galois if $b \ne 0$. 
\end{proof} 

\
\begin{center} {\bf Acknowledgements} \end{center} 
The authors are grateful to Professor Hisao Yoshihara for helpful discussions. 
In particular, Yoshihara informed the first author an idea of the proof of Theorem \ref{HyperplaneSection} when $p=0$, $s(X) \le 0$ and $H$ is general.  
The first author was supported by Research Fellowships of the Japan Society for the Promotion of Science for Young Scientists.

\end{document}